\documentclass[12pt,a4paper]{amsart}

\usepackage[utf8]{inputenc}
\usepackage[american]{babel}
\usepackage{amsmath,amssymb,amsthm}
\usepackage[shortlabels]{enumitem}
\usepackage{graphicx} 
\usepackage{braket}
\usepackage[left=3cm,right=3cm]{geometry}
\usepackage[hidelinks]{hyperref}

\newcommand{\Cay}{\operatorname{Cay}}

\newcommand{\abs}[1]{\lvert#1\rvert}
\newcommand{\norm}[1]{\lVert#1\rVert}

\newtheorem{theorem}{Theorem}[section]
\newtheorem{proposition}[theorem]{Proposition}
\newtheorem{lemma}[theorem]{Lemma}

\theoremstyle{remark}
\newtheorem{remark}[theorem]{Remark}

\title{Rigidity in the Lichnerowicz inequality for word-length semigroups}
\author{Melchior Wirth}
\address{Faculty of Mathematics and Computer Science, Leipzig University, Neues Augusteum, Augustusplatz 10, 04109 Leipzig, Germany}
\email{melchior.wirth@uni-leipzig.de}

\begin{document}

\begin{abstract}
    We give an operator-algebraic characterization of Cayley graphs that are partial cubes, that is, embed isometrically into a hypercube. More precisely, noncommutative versions of the Lichnerowicz inequality control the spectral gap of the generator of a quantum Markov semigroup by the constant in the complete $\Gamma_2$-criterion of Junge and Zeng or the constant in the geometric Ricci criterion of Junge, LaRacuente and Li. We show that for word-length semigroups on group von Neumann algebras, equality in either of these inequalities occurs exactly when the Cayley graph is a partial cube.
\end{abstract}

\maketitle

\section{Introduction}

In Riemannian geometry, the Lichnerowicz inequality \cite{Lic58} gives a lower bound for the spectral gap of the Laplace--Beltrami operator in terms of a lower bound on the Ricci curvature. In recent years, a number of approaches have been developed to extend the notion of lower Ricci curvature bounds beyond Riemannian manifolds, including the setting of metric measure spaces, discrete graphs and general (quantum) Markov semigroups. While not all of these cases support a notion of dimension, in each case a dimension-free version of the Lichnerowicz inequality is known \cite{BE85,LV09,JZ15,LJL25}.

A natural question is to ask when the Lichnerowicz inequality is saturated. For the dimension-dependent inequality for Riemannian manifolds, Obata's theorem \cite{Oba62} states that the equality cases are exactly the round spheres. More generally, for RCD spaces, the equality cases are exactly the spherical suspensions by a result of Ketterer \cite{Ket15}. Recently, versions of Obata's rigidity have also been obtained for discrete graphs and quantum Markov semigroups induced by free difference quotients \cite{Die26,DLZ26a,DLZ26b}.

In this article, we characterize the equality cases in the Lichnerowicz inequality for word-length semigroups on group von Neumann algebras that arise from a generating set of the group. In this setting, there are several (potentially inequivalent) notions of lower Ricci curvature bounds, including the complete $\Gamma_2$-criterion of Junge and Zeng \cite{JZ15} and the stronger geometric Ricci condition of Li, Junge and LaRacuente \cite{LJL25}, which all imply a Lichnerowicz inequality for the spectral gap of the generator.

We show that for both noncommutative Ricci curvature notions mentioned in the previous paragraph, the equality cases in the Lichnerowicz inequality are the same and are characterized by the fact that the Cayley graph embeds isometrically into a hypercube. Graphs with this property are called partial cubes.

This equivalence is particularly interesting as these noncommutative Ricci curvature notions are defined by abstract analytic criteria concerning the Markov semigroups and their generators, and so far, no direct geometric implications for the geometry of the group had been known. This opens up the field to explore connections between analytic properties of word-length semigroups and the geometry of the underlying Cayley graph in general.

Moreover, it is an open question whether the complete $\Gamma_2$-criterion implies the geometric Ricci condition for word-length semigroups. One consequence of our result is that this implication is true at least in the extreme case when the Lichnerowicz inequality is saturated.

For each curvature notion discussed above, both implications -- that word-length semigroups of Cayley graphs that are partial cubes saturate the Lichnerowicz inequality, and that all equality cases have Cayley graphs that are partial cubes -- are new. While it is known that all word-length semigroups satisfy the geometric Ricci condition with constant $0$, obtaining examples with strictly positive constant has proven difficult so far. One of the starting points of this investigation was the observation that most known examples with strictly positive curvature in the sense of the geometric Ricci condition can be verified using hypercube transference and saturate the Lichnerowicz inequality. As our main result and its proof show, all equality cases in the Lichnerowicz inequality arise this way.

Cayley graphs that are partial cubes include the Cayley graphs of Coxeter groups with their defining generating set. It had been conjectured that all finite Cayley graphs that are partial cubes are Cayley graphs of Coxeter groups, but this conjecture has been recently refuted by Knauer and Marc (see \cite[Chapter 3]{Jae25}). The fact that word-length functions of Coxeter systems saturate the Lichnerowicz inequality has been known in a number of cases, including free groups and products of cyclic groups of order $2$, but is new in full generality. While all previous proofs relied either on the (free) product structure or approximation properties such as amenability, our result shows that neither of these properties is necessary.

The proof of our main result rests on the following steps: We first show that the geometric Ricci condition is stable under wreath products. Then we prove that the partial cube property of the Cayley graph induces an embedding of the group into a wreath product of $C_2$ with the group that intertwines the semigroups. Since it is known that the word-length semigroup on $C_2$ satisfies the geometric Ricci condition with constant $1$, this implies that the same is true for any word-length semigroup for which the Cayley graph is a partial cube.

For the converse implication, the optimality condition implies that the length function changes parity along every edge. It follows that the graph is bipartite. Since Cayley graphs are connected and we assume the length function to be of negative type, it follows from known results in graph theory that the Cayley graph is a partial cube.

This article is structured as follows. In Section~\ref{sec:prelim} we gather the relevant definitions and results concerning conditionally negative functions on groups, word-length semigroups and Ricci curvature notions. In Section~\ref{sec:main} we show that for amenable groups, the $\Gamma_2$-criterion and the complete $\Gamma_2$-criterion coincide, before we state and prove the main result of this article (Theorem~\ref{thm:main}), in which we characterize the extreme cases in the noncommutative Lichnerowicz inequality as the word-length semigroups arising from Cayley graphs that are partial cubes.

\subsection*{Acknowledgments} The author is grateful to Marius Junge for inspiring discussions on noncommutative Ricci curvature on various occasions. He gratefully acknowledges funding by the Deutsche Forschungsgemeinschaft (DFG, German Research Foundation) -- Project number 565405322.

The author used ChatGPT-6 Astra to explore the connections to graph theory and in particular the literature on partial cubes. The text of this article was written by the author and he takes full responsibility for its content.

\section{Preliminaries}\label{sec:prelim}

In this section we give a brief overview of some background material. For more details, we refer to \cite{BdlHV08} for conditionally negative definite functions on groups, to \cite{JZ15,LJL25} for noncommutative Ricci curvature notions and to \cite{Djo73} for partial cubes.

Let $G$ be a (discrete) group. A function $\varphi\colon G\to \mathbb C$ is called \emph{positive definite} if
\begin{equation*}
    \sum_{g,h\in G}\overline{c_g}c_h\varphi(g^{-1}h)\geq 0
\end{equation*}
whenever $(c_g)_{g\in G}$ is a finitely supported family in $\mathbb C$.

A function $\ell\colon G\to [0,\infty)$ is called \emph{conditionally negative definite (cnd)} if $\ell(e)=0$, $\ell(g^{-1})=\ell(g)$ for all $g\in G$ and
\begin{equation*}
    \sum_{g,h}\overline{c_g}c_h\ell(g^{-1}h)\leq 0
\end{equation*}
whenever $(c_g)_{g\in G}$ is a finitely supported family in $\mathbb C$ such that $\sum_{g\in G}c_g=0$. By Schönberg's theorem, the function $\ell$ is cnd if and only if $e^{-t\ell}$ is positive definite for all $t\geq 0$.

If $\ell\colon G\to [0,\infty)$ is a cnd function and $t\geq 0$, then there exists a unique normal unital completely positive map $T_t^\ell$ on the group von Neumann algebra $L(G)$ such that $T_t^\ell(\lambda(g))=e^{-t\ell(g)}\lambda(g)$ for all $g\in G$. Moreover, the family $(T_t^\ell)_{t\geq 0}$ is an operator semigroup, continuous in the point-weak$^\ast$ topology, and
\begin{equation*}\label{eq:tracial_sym}
    \tau(T_t^\ell(x)^\ast y)=\tau(x^\ast T_t^\ell(y))\tag{$\ast$}
\end{equation*}
for all $x,y\in L(G)$ and $t\geq 0$. Here $\tau$ denotes the canonical tracial state on $L(G)$. In general, if $M$ is a finite von Neumann algebra with faithful normal tracial state $\tau$, a point-weak$^\ast$ continuous operator semigroup of normal unital completely positive maps on $M$ is called a \emph{quantum Markov semigroup}. It is called \emph{$\tau$-symmetric} if \eqref{eq:tracial_sym} holds.

In this article we will focus on cnd functions that arise as word-length functions of generating sets. If $S\subset G\setminus\{e\}$ is a symmetric generating set, the \emph{word-length function} is defined by
\begin{equation*}
    \ell_S(g)= \min\{n\in \mathbb N_0\mid \exists s_1,\dots,s_n\in S:g=s_1\dots s_n\}
\end{equation*}
for $g\neq e$ and $\ell_S(e)=0$. The \emph{combinatorial graph distance} on the Cayley graph $\Cay(G,S)$ is given by $d(g,h)=\ell_S(g^{-1}h)$ for $g,h\in G$.

If $I$ is an index set and $i\in I$, let $s_i$ denote the image of the generator of the cyclic group $C_2$ in the $i$-th copy of $C_2$ inside $C_2^{(I)}=\bigoplus_{i\in I}C_2$. The set $S_I=\{s_i\mid i\in I\}$ is a symmetric generating set of $C_2^{(I)}$. A \emph{hypercube} is a Cayley graph $\Cay(C_2^{(I)},S_I)$ for some (possibly infinite) index set $I$. A graph is a \emph{partial cube} if it embeds isometrically with respect to the combinatorial graph distances into a hypercube. 

In general, word-length functions are not necessarily cnd, but they are cnd if the Cayley graph $\Cay(G,S)$ is a partial cube, as we will see below. In this case, the word-length semigroup $(T_t^{\ell_S})_{t\geq 0}$ on $L(G)$ is a quantum Markov semigroup, and one can study various notions of Ricci curvature bounds. For the Bakry--Émery or $\Gamma_2$-criterion, there is a direct formulation in terms of the length function following Junge and Zeng \cite{JZ15}.

Let us introduce the criterion more generally for an arbitrary cnd function $\ell$ on $G$. Define
\begin{equation*}
    K\colon G\times G\to \mathbb R,\,K(g,h)=\frac 1 2(\ell(g)+\ell(h)-\ell(g^{-1}h)).
\end{equation*}
The semigroup $(T_t^\ell)_{t\geq 0}$ satisfies the \emph{$\Gamma_2$-criterion with constant $\kappa$} if $\sum_{g,h}(K(g,h)^2-\kappa K(g,h))\overline{c_g}c_h\lambda(g^{-1}h)$ is positive in $L(G)$ for every finitely supported family $(c_g)_{g\in G}$ in $\mathbb C$. It satisfies the \emph{complete $\Gamma_2$-criterion with constant $\kappa$} if $\sum_{g,h}(K(g,h)^2-\kappa K(g,h)) C_g^\ast C_h\otimes\lambda(g^{-1}h)$ is positive in $M_n(\mathbb C)\otimes L(G)$ for every $n\in\mathbb N$ and every finitely supported family $(C_g)_{g\in G}$ in $M_n(\mathbb C)$.

The complete $\Gamma_2$-criterion has the following simpler formulation, which is the original definition by Junge and Zeng \cite{JZ15}. Since the equivalence is not explicitly stated there, we include its simple proof. We write $\circ $ for the Hadamard product of (possibly infinite) matrices and $A\geq 0$ for an infinite matrix means that the blocks $[A(j,k)]_{j,k\in F}$ are positive for every finite index set $F$.

\begin{lemma}\label{lem:char_CBE}
    The semigroup $(T_t^\ell)_{t\geq 0}$ satisfies the complete $\Gamma_2$-criterion with constant $\kappa$ if and only if $K\circ K\geq \kappa K$.
\end{lemma}
\begin{proof}
    First assume that the complete $\Gamma_2$-criterion is satisfied and let $A=K\circ K-\kappa K$. Let $F=\{g_1,\dots,g_n\}$, let $C_{g_j}=\ket 1\bra j$ for $1\leq j\leq n$ and let $C_g=0$ for $g\notin F$. The hypothesis implies $[A(g_j,g_k)\lambda(g_j^{-1}g_k)]_{j,k=1}^n\geq 0$. With the unitary matrix $U=\mathrm{diag}(\lambda(g_1),\dots,\lambda(g_n))$ one gets
    \begin{equation*}
        0\leq [A(g_j,g_k)\lambda(g_j^{-1}g_k)]_{j,k=1}^n=U^\ast ([A(g_j,g_k)]_{j,k=1}^n\otimes 1)U.
    \end{equation*}
    Thus $A=K\circ K-\kappa K\geq 0$.

    Conversely, if $K\circ K-\kappa K\geq 0$, $F\subset G$ is finite and $(C_g)_{g\in F}$ is a family in $M_n(\mathbb C)$, let $B\in M_F(\mathbb C)$ such that $(B^\ast B)(g,h)=(K\circ K-\kappa K)(g,h)$ for $g,h\in F$. We have
    \begin{equation*}
        \sum_{g,h}(K(g,h)^2-\kappa K(g,h)) C_g^\ast C_h\otimes\lambda(g^{-1}h)=\sum_{s\in F}\left\lvert B_{s,g}C_g\otimes \lambda(g)\right\rvert^2\geq 0.\qedhere
    \end{equation*}
\end{proof}

For the formulation of the geometric Ricci condition, we need some more operator-algebraic machinery. Every $\tau$-symmetric quantum Markov semigroup $(T_t)_{t\geq 0}$ on a finite von Neumann algebra $M$ extends uniquely to a strongly continuous symmetric contraction semigroup on $L^2(M,\tau)$. Thus there exists a unique positive self-adjoint operator $L$ in $L^2(M,\tau)$ such that $T_t=e^{-tL}$.

A \emph{derivation triple} for $(T_t)_{t\geq 0}$ is a triple $(M\subset \hat M,\hat \tau,\delta)$ consisting of a unital inclusion $M\subset \hat M$ of von Neumann algebras, a faithful normal tracial state $\hat\tau$ on $\hat M$ that extends $\tau$ and a real closable derivation $\delta\colon \mathrm{dom}(\delta)\subset L^2(M,\tau)\cap M\to L^2(\hat M,\hat \tau)$ such that $\delta^\ast\bar\delta=L$.

Write $\mathcal H_\delta$ for the closed $M$-bimodule generated by $\operatorname{ran}\delta$. The semigroup $(T_t)_{t\geq 0}$ is said to satisfy the \emph{geometric Ricci condition with constant $\kappa$} if there exists a derivation triple $(M\subset\hat M,\hat\tau,\delta)$, a $\hat\tau$-symmetric quantum Markov semigroup $(e^{-t\hat L})_{t\geq 0}$ on $\hat M$ and a bounded $M$-bimodule map $R\colon \mathcal H_\delta\to L^2(\hat M,\hat\tau)$ such that
\begin{equation*}
    \bar\delta T_t=e^{-t(\hat L+R)}\bar\delta
\end{equation*}
on $\mathrm{dom}(\bar\delta)$ for $t>0$ and $\langle\xi,R\xi\rangle\geq \kappa\norm{\xi}^2$ for $\xi\in\mathcal H_\delta$.

The following permanence property is immediate from the definition.

\begin{lemma}\label{lem:GRic_intertwining}
    For $j\in \{1,2\}$ let $M_j$ be a finite von Neumann algebra, $\tau_j$ a faithful normal tracial state on $M_j$ and $(T_t^j)_{t\geq 0}$ a $\tau_j$-symmetric quantum Markov semigroup on $M_j$. If $\pi\colon M_1\to M_2$ is a trace-preserving $\ast$-homomorphism such that $\pi\circ T_t^1=T_t^2\circ\pi$ for $t\geq 0$ and $(T_t^2)_{t\geq 0}$ satisfies the geometric Ricci condition with constant $\kappa$, so does $(T_t^1)_{t\geq 0}$.
\end{lemma}

\section{Main result}\label{sec:main}

In this section we state and prove the main result of this article, the characterization of the equality cases in the Lichnerowicz inequality for word-length semigroups as exactly the semigroups induced by Cayley graphs that are partial cubes. We also show that the $\Gamma_2$-criterion and the complete $\Gamma_2$-criterion coincide for amenable groups and arbitrary cnd functions.

\begin{proposition}\label{prop:BE_amenable}
    If $G$ is an amenable group, $\ell\colon G\to [0,\infty)$ is a cnd function, and $\kappa\geq 0$, then the following assertions are equivalent:
    \begin{enumerate}[(i)]
        \item $(T_t^\ell)_{t\geq 0}$ satisfies the $\Gamma_2$-criterion with constant $\kappa$
        \item $(T_t^\ell)_{t\geq 0}$ satisfies the complete $\Gamma_2$-criterion with constant $\kappa$
    \end{enumerate}
\end{proposition}
\begin{proof}
    The implication (ii)$\implies$(i) is clear. For the implication (i)$\implies$(ii), let $(c_g)_{g\in G}$ be a finitely supported family in $\mathbb C$, $x=\sum_g c_g \lambda_g$ and $A=K\circ K-\kappa K$. The $\Gamma_2$-criterion with constant $\kappa$ implies 
    \begin{equation*}
        \sum_{g,h}\overline{c_g}c_h A(g,h)\lambda_{g^{-1}h}=\Gamma_2(x)-\kappa \Gamma(x)\geq 0.
    \end{equation*}
    If $(F_i)_{i\in I}$ is a F{\o}lner net, then
    \begin{align*}
        \frac 1{\abs{F_i}}\langle 1_{F_i},(\Gamma_2(x)-\kappa \Gamma(x))1_{F_i}\rangle=\sum_{g,h}\overline{c_g}c_h A(g,h)\frac{\abs{F_i\cap g^{-1}h F_i}}{\abs{F_i}}\overset{i}{\to}\sum_{g,h}\overline{c_g}c_h A(g,h).
    \end{align*}
    Thus $A\geq 0$, which is equivalent to the complete $\Gamma_2$-criterion with constant $\kappa$ by Lemma~\ref{lem:char_CBE}.
\end{proof}

Let us record some simple permanence properties of cnd functions.

\begin{lemma}\label{lem:permanence_cnd}
    \begin{enumerate}[(a)]
        \item If $G,\Gamma$ are groups, $\alpha\colon G\to\mathrm{Aut}(\Gamma)$ is a group homomorphism and $\ell\colon \Gamma\to [0,\infty)$ is a cnd function such that $\ell\circ\alpha_g=\ell$ for all $g\in G$, then the function 
        \begin{equation*}
            \ell^\rtimes\colon \Gamma\rtimes_\alpha G\to [0,\infty),\,\ell^\rtimes(\gamma g)=\ell(\gamma)\qquad\gamma\in\Gamma,g\in G,
        \end{equation*}
        is cnd.

        \item If $(\Gamma_i)_{i\in I}$ is a family of groups and $\ell_i\colon \Gamma_i\to [0,\infty)$ is a cnd function, then the function 
        \begin{equation*}
            \ell\colon \bigoplus_{i\in I}\Gamma_i\to [0,\infty),\,\ell((\gamma_i)_{i\in I})=\sum_{i\in I}\ell_i(\gamma_i),\qquad(\gamma_i)_{i\in I}\in\bigoplus_{i\in I}\Gamma_i,
        \end{equation*}
        is cnd.
    \end{enumerate}
\end{lemma}
\begin{proof}
    \begin{enumerate}[(a)]
        \item The relations of the semidirect product and the invariance of $\ell$ imply
        \begin{equation*}
            \ell^\rtimes(g\gamma h)=\ell^\rtimes(\alpha_g(\gamma)gh)=\ell^\rtimes(\alpha_g(\gamma))=\ell(\gamma)
        \end{equation*}
        for all $g,h\in G$, $\gamma\in \Gamma$.

        Now, if $(c_{\gamma,g})_{\gamma\in \Gamma,g\in G}$ is a finitely supported family in $\mathbb C$ such that $\sum_{\gamma,g}c_{\gamma,g}=0$, let $c_\gamma=\sum_{g\in G}c_{\gamma,g}$ for $\gamma\in \Gamma$. We have
        \begin{equation*}
            \sum_{\gamma,\delta,g,h}\overline{c_{\gamma,g}}c_{\delta,h}\ell^\rtimes((\gamma g)^{-1}(\delta h))=\sum_{\gamma,\delta,g,h}\overline{c_{\gamma,g}}c_{\delta,h}\ell(\gamma^{-1}\delta)=\sum_{\gamma,\delta}\overline{c_\gamma}c_\delta\ell(\gamma^{-1}\delta)\leq 0.
        \end{equation*}
        \item Let $\Gamma=\bigoplus_{i\in I}\Gamma_i$. If $\gamma^1,\dots,\gamma^n\in \Gamma$ and $c_1,\dots,c_n\in\mathbb C$ such that $\sum_{j=1}^n c_j=0$, then
        \begin{equation*}
            \sum_{j,k=1}^n\overline{c_j}c_k \ell((\gamma^j)^{-1}\gamma^k)=\sum_{i\in I}\sum_{j,k=1}^n\overline{c_j}c_k\ell_i((\gamma^j_i)^{-1}\gamma^k_i)\leq 0.\qedhere
        \end{equation*}
    \end{enumerate}
\end{proof}

If $\Gamma$, $G$ are groups, $I$ is an index set and $G$ acts on $I$ by permutations, we write $\Gamma \wr_I G$ for the \emph{restricted wreath product}, which is the crossed product $(\bigoplus_{i\in I}\Gamma)\rtimes G$ for the permutation action of $G$.

\begin{lemma}\label{lem:GRic_wreath}
    If $\Gamma$, $G$ are groups, $I$ is an index set, $G$ acts on $I$ by permutations and $\ell\colon \Gamma\to [0,\infty)$ is a cnd function, then the function 
    \begin{equation*}
        \ell^\wr\colon \Gamma\wr_I G\to [0,\infty),\,\ell^\wr((\gamma_i)_{i\in I},g)=\sum_{i\in I}\ell(\gamma_i),\qquad(\gamma_i)_{i\in I}\in \bigoplus_{i\in I}\Gamma,g\in G,
    \end{equation*}
    is cnd.

    Moreover, if $(T_t^\ell)_{t\geq 0}$ satisfies the geometric Ricci condition with constant $\kappa$, so does $(T_t^{\ell^\wr})_{t\geq 0}$.
\end{lemma}
\begin{proof}
    The fact that $\ell^\wr$ is cnd follows from Lemma~\ref{lem:permanence_cnd}. Let $(L(\Gamma)\subset \hat M,\hat\tau,\delta)$ be a derivation triple for $(T_t^\ell)_{t\geq 0}$, let $(e^{-t\hat L})_{t\geq 0}$ be a $\hat\tau$-symmetric quantum Markov semigroup on $\hat M$ and let $R\colon \mathcal H_\delta\to L^2(\hat M,\hat \tau)$ be a bounded $L(\Gamma)$-bimodule map such that $\langle \xi,R\xi\rangle\geq \kappa\norm{\xi}^2$ for $\xi\in\mathcal H_\delta$ and $\delta T_t^\ell=e^{-t(\hat L+R)}\delta$ on $\mathrm{dom}(\delta)$.

    Consider the tracial tensor product $\bigotimes_{i\in I}(\hat M,\hat \tau)$, on which $G$ acts by coordinate permutations.  We have an inclusion
    \begin{equation*}
        L(\Gamma\wr_I G)\cong L\bigg(\bigoplus_{i\in I}\Gamma\bigg)\rtimes G\cong \bigotimes_{i\in I}(L(\Gamma),\tau)\rtimes G\subset \bigotimes_{i\in I}(\hat M,\hat \tau)\rtimes G.
    \end{equation*}
    The $L^2$-space corresponding to $L(\Gamma\wr_I G)$ is then given by $\left(\bigotimes_{i\in I}L^2(L(\Gamma),\tau)\right)\otimes\ell^2(G)$, where the infinite tensor product is the closed linear hull of all elementary tensors $\bigotimes_{i\in I}\xi_i$ for which $\xi_i=1$ for all but finitely indices $i$. Furthermore, since $G$ acts by trace-preserving $\ast$-automorphisms on $\bigotimes_{i\in I}(\hat M,\hat \tau)$, the von Neumann algebra on the right side has a faithful normal tracial state $\mathcal T$ that extends the canonical trace on $L(\Gamma\wr_I G)$.

    Let $\mathrm{dom}(\partial)$ be the linear hull of the elementary tensors $\bigotimes_{i\in I}\xi_i\otimes x\delta_e$ with $\xi_i\in \mathrm{dom}(\delta)$, $\xi_i=1$ for all but finitely many $i\in I$ and $x\in L(G)$. Define
    \begin{equation*}
        \partial\colon \mathrm{dom}(\partial)\to \bigotimes_{i\in I}L^2(\hat M,\hat\tau)\otimes \ell^2(G),\,\partial\left(\bigotimes_{i\in I}\xi_i\otimes x\delta_e\right)=\sum_{i\in I}\bigotimes_{j\neq i}\xi_j\otimes \delta(\xi_i)\otimes x\delta_e,
    \end{equation*}
    where the sum on the right side is finite by assumption. Note that $\mathcal H_\partial$ is the closed linear subspace of $\bigotimes_{i\in I}L^2(\hat M,\hat \tau)\otimes\ell^2(G)$ generated by the elementary tensors for which exactly one factor is in $\mathcal H_\delta$ and the other factors are in $L^2(L(\Gamma),\tau)$.

    Moreover, define $\mathcal L$ by $e^{-t\mathcal L}=\bigotimes_{i\in I}e^{-t\hat L}\otimes \mathrm{id}_{\ell^2(G)}$ and $\mathcal R(\bigotimes_{i\in I}\xi_i\otimes \eta)=\bigotimes_{i\neq j}\xi_i\otimes R\xi_j\otimes\eta$ if $\xi_i\in L^2(L(\Gamma),\tau)$ for $i\neq j$, $\xi_j\in \mathcal H_\delta$ and $\eta\in \ell^2(G)$. A direct computation shows that $(L(\Gamma\wr_I G)\subset \bigotimes_{i\in I}(\hat M,\hat\tau)\rtimes G,\mathcal T,\partial)$ is a derivation triple for $(T_t^{\ell^\wr})_{t\geq 0}$ and that $\mathcal L$, $\mathcal R$ witness the geometric Ricci criterion with constant $\kappa$.
\end{proof}

We will use the following relation, which plays an important role in the study of partial cubes. The \emph{Winkler relation} \cite{Win84} on the edges of $\Cay(G,S)$ is defined by $(g,h)\Theta(g^\prime,h^\prime)$ if
\begin{equation*}
    d(g,g^\prime)+d(h,h^\prime)\neq d(g,h^\prime)+d(h,g^\prime)
\end{equation*}
for $g,h,g^\prime,h^\prime\in G$ with $g^{-1}h,(g^\prime)^{-1} h^\prime\in S$. Clearly, the Winkler relation is invariant under isometric embeddings and the left translation action of $G$ on the edges of $\Cay(G,S)$.

For $\Cay(C_2^{(I)},S_I)$ one has $(g,h)\Theta(g^\prime,h^\prime)$ if and only if $g^{-1}h=(g^\prime)^{-1} h^\prime$. To see this, realize the vertices of $\Cay(C_2^{(I)},S_I)$ as elements of $\ell^2(I)$ with coordinates in $\{0,1\}$ and note that $d(v,w)=\norm{v-w}_2^2$ for this realization. Then one has $v-w=\pm e_i$, $v^\prime-w^\prime=\pm e_j$ for neighboring vertices and
\begin{equation*}
    d(v,v^\prime)+d(w,w^\prime)- d(v,w^\prime)-d(w,v^\prime)=-2\langle v-w,v^\prime-w^\prime\rangle,
\end{equation*}
which is non-zero exactly if $i=j$.

\begin{theorem}\label{thm:main}
    For a group $G$ with symmetric generating set $S\subset G\setminus\{e\}$, the following assertions are equivalent:
    \begin{enumerate}[(i)]
        \item $\Cay(G,S)$ is a partial cube
        \item $\ell_S$ is cnd and $(T_t^{\ell_S})_{t\geq 0}$ satisfies the geometric Ricci condition with constant $1$
        \item $\ell_S$ is cnd and $(T_t^{\ell_S})_{t\geq 0}$ satisfies the complete $\Gamma_2$-criterion with constant $1$
    \end{enumerate}
    If $G$ is amenable, these assertions are also equivalent to
    \begin{enumerate}[(i)]
        \setcounter{enumi}{3}
        \item $\ell_S$ is cnd and $(T_t^{\ell_S})_{t\geq 0}$ satisfies the $\Gamma_2$-criterion with constant $1$
    \end{enumerate}
\end{theorem}
\begin{proof}
    (i)$\implies$(ii): Assume that $\phi\colon\Cay(G,S)\to\Cay(\bigoplus_{i\in I}C_2,S_I)$ is an isometric hypercube embedding and let $\ell$ denote the word-length function on $C_2$. As discussed above, if we realize $\Cay(\bigoplus_{i\in I}C_2,S_I)$ inside $\ell^2(I)$, we get that
    \begin{equation*}
        \ell_S(g^{-1}h)=d(g,h)=\norm{\phi(g)-\phi(h)}_2^2
    \end{equation*}
    for $g,h\in G$. This immediately implies that $\ell_S$ is cnd.
    
    To see that $(T_t^{\ell_S})_{t\geq 0}$ satisfies the geometric Ricci condition with constant $1$, we will construct a permutation action of $G$ on $I$ and an injective group homomorphism $\alpha\colon G\to C_2\wr_I G$ such that $\ell_S=\ell^\wr\circ \alpha$.

    If $h_1,h_2\in G$ such that $h_1^{-1}h_2\in S$, then $\phi(h_1)^{-1}\phi(h_2)\in S_I$. Let $i(h_1,h_2)$ be the unique index $i\in I$ such that $\phi(h_1)^{-1}\phi(h_2)=s_i$. As discussed above, $i(h_1,h_2)=i(h_1^\prime,h_2^\prime)$ if and only if $(h_1,h_2)\Theta(h_1^\prime,h_2^\prime)$ and the Winkler relation is invariant under the left translation action on the edges. Thus we can define a permutation action of $G$ on $I$ by $\sigma_g(i(h_1,h_2))=i(gh_1,gh_2)$ and $\sigma_g(i)=i$ if there is no edge $(h_1,h_2)$ such that $i(h_1,h_2)=i$.
    
    To define the embedding $\alpha\colon G\to C_2\wr_I G$, let $b(g)=\sum_{i:\phi(e)_i\neq \phi(g)_i}s_i$ and $\alpha(g)=b(g)g$ for $g\in G$. If $g,h\in G$ and $i\in I$, then $\phi(e)_i\neq \phi(gh)_i$ if and only if either $\phi(e)_i\neq \phi(g)_i$ or $\phi(g)_i\neq \phi(gh)_i$. Hence
    \begin{equation*}
        b(gh)=b(g)+\sum_{i: \phi(g)_i\neq \phi(gh)_i}s_i=b(g)+\sigma_g b(h).
    \end{equation*}
    In other words, $b$ is a $1$-cocycle relative to $\sigma$. Therefore, $\alpha$ is a group homomorphism. It is injective since $g\neq e$ implies $b(g)g\neq e$ in $C_2\wr_I G$. We get a trace-preserving unital $\ast$-homomorphism $\pi\colon L(G)\to L(C_2\wr_I G)$ defined by $\pi(\lambda(g))=\lambda(\alpha(g))$ for $g\in G$.

    By \cite[Theorem 5.3 (ii)]{BGJ23} (see also Example 5.4), the word-length semigroup $(T_t^\ell)_{t\geq 0}$ satisfies the geometric Ricci condition with constant $1$. By Lemma~\ref{lem:GRic_wreath}, the same is then also true for $(T_t^{\ell^\wr})_{t\geq 0}$. Finally, this continues to hold for $(T_t^{\ell_S})_{t\geq 0}$ by Lemma~\ref{lem:GRic_intertwining}.

    (ii)$\implies$(iii): This follows from \cite[Proposition 2.7]{BGJ23}.

    (iii)$\implies$(i): Let $A=(K(g,h)^2-K(g,h))_{g,h\in G}$. By (iii), $A\geq 0$. If $s\in S$, then $A(s,s)=0$, which implies $A(g,s)=A(s,g)=0$ and thus $K(g,s)=K(s,g)\in \{0,1\}$ for all $g\in G$. Hence
    \begin{equation*}
        \ell_S(g)-\ell_S(gs)=\ell_S(g^{-1})-\ell_S((g^{-1})^{-1}s)=2 K(g^{-1},s)-\ell_S(s)\in \{-1,1\}.
    \end{equation*}
    If we define $G_j=\{g\in G\mid \ell_S(g)\equiv j\pmod{2}\}$ for $j\in \{0,1\}$, it follows that there are no edges in $\Cay(G,S)$ with both endpoints in $G_0$ or both endpoints in $G_1$. Therefore, $\Cay(G,S)$ is bipartite.
    
    Moreover, $\Cay(G,S)$ is connected, and since $\ell_S$ is cnd, there exists a real Hilbert space $H$ and a map $b\colon G\to H$ such that $d(g,h)=\norm{b(g)-b(h)}^2$ for $g,h\in G$ by \cite[Theorem C.2.3]{BdlHV08}. Now it follows from \cite[Section 3]{RW86} (see also \cite[Lemma 1, Theorems 2, 3]{Tan23} for the case of infinite graphs) that $\Cay(G,S)$ is a partial cube.

    If $G$ is amenable, the equivalence (iii)$\iff$(iv) follows from Proposition~\ref{prop:BE_amenable}.
\end{proof}

\begin{remark}
    In terms of the Gromov form $K$, the equivalence between (i) and (iii) can be summarized as follows:
    \begin{equation*}
        \Cay(G,S)\text{ is a partial cube}\text{ if and only if }K\circ K\geq K\geq 0.
    \end{equation*}
    This is a directly accessible criterion that does not require operator-algebraic machinery.
\end{remark}

\begin{remark}
    As one sees from the proof of (iii)$\implies$(i), if $\ell_S$ is cnd, then $\Cay(G,S)$ is a partial cube if and only if $\ell_S$ changes parity along every edge.
\end{remark}

\bibliographystyle{alphaurl}
\bibliography{ref}

\end{document}